\documentclass[12pt]{article}
\usepackage[utf8]{inputenc}
\usepackage{amssymb,amsmath,amsfonts,eucal,mathrsfs,amsthm} 
\usepackage[allcolors=blue]{hyperref}
\usepackage{empheq}
\newtheorem{theorem}{Theorem}

\newtheorem{example}[theorem]{Example}

\theoremstyle{definition}
\newcommand{\R}{\mathbb{R}}
\newcommand{\Q}{\mathbb{Q}}
\newcommand{\Sf}{\mathbb{S}}

\newcommand{\Hy}{\mathbb{H}}
\newcommand{\spa}{\mbox{span}}

\def\<{{\langle}}
\def\>{{\rangle}}

\def\a{\alpha}

\def\be{\begin{equation} }
\def\ee{\end{equation} }

\def\qed{\ifhmode\unskip\nobreak\fi\ifmmode\ifinner
\else\hskip5 pt \fi\fi\hbox{\hskip5 pt \vrule width4 pt
height6 pt  depth1.5 pt \hskip 1pt }}
\begin{document}

\title{On constant curvature submanifolds\\ of space forms}
\author{M. Dajczer, C.-R. Onti and Th. Vlachos}
\date{}
\maketitle

\begin{abstract} We prove a converse to well-known results by 
E. Cartan and J. D. Moore.  Let $f\colon M^n_c\to\Q^{n+p}_{\tilde c}$ 
be an isometric immersion of a Riemannian manifold with constant 
sectional curvature $c$ into a space form of curvature $\tilde c$, 
and free of weak-umbilic points if $c>\tilde{c}$. We show that the 
substantial codimension of $f$ is  $p=n-1$ if, as shown by Cartan 
and Moore, the first normal bundle possesses the lowest possible 
rank $n-1$. These submanifolds are of a class that has been extensively 
studied due to their many properties. For instance, they are 
holonomic and admit B\"{a}cklund and  Ribaucour transformations. 
\end{abstract}

E. Cartan \cite{ca} in 1919 initiated the systematic study 
of isometric immersions $f\colon M^n_c\to\Q^{n+p}_{\tilde c}$ 
of an $n$-dimensional connected Riemannian manifold $M^n_c$ of 
constant sectional curvature $c$ into a simply connected space 
form. Thus $\Q^m_{\tilde c}$ 
denotes the Euclidean space $\R^m$, the Euclidean sphere 
$\Sf_{\tilde c}^m$ or the hyperbolic space  $\Hy_{\tilde c}^m$ 
according to whether  $\tilde c=0,\tilde c>0$ 
or $\tilde c<0$, respectively. 

Cartan proved for $n\geq 3$ and $c<\tilde c$, that the 
codimension of $f$ satisfies $p\geq n-1$, and if $p=n-1$  
the normal bundle has to be flat. The dual case of isometric 
immersions $f\colon M^n_c\to\Q^{n+p}_{\tilde c}$, $n\geq 3$ 
and $c>\tilde c$, has been considered by J. D. Moore \cite{mo}. 
Under the additional assumption that $f$ is free of weak-umbilic 
points, he again obtained that $p\geq n-1$, and if $p=n-1$
that the normal bundle is flat . We recall 
that $x\in M_c^n$ is called a \emph{weak-umbilic} for $f$ if 
there exists a unit normal vector $\zeta\in N_fM(x)$ such that 
the corresponding shape operator is $A_\zeta=\sqrt{c-\tilde c}\ I$.
\vspace{1ex}

The way to argue that $p$ satisfies the lower bound given above 
is to show that the dimension of the first normal space has to 
satisfy $\dim N_1^f\geq n-1$ at any point. Recall that the 
\emph{first normal space} $N_1^f(x)$ of $f$ at $x\in M^n$ 
is the subspace of the normal space $N_fM(x)$ spanned by the 
image of its second fundamental form 
$\a_f\colon TM\times TM\to N_fM$ at $x\in M^n$, that is,
$$
N_1^f(x)=\spa\{\a_f(X,Y):X,Y\in T_xM\}.
$$

The purpose of this paper is to give a converse for isometric 
immersions $f\colon M^n_c\to\Q^{n+p}_{\tilde c}$ of the 
aforementioned results by Cartan and Moore,  that goes as 
follows:  If $n\geq 2$ and assuming that $\dim N_1^f=n-1$ at any point, 
we show that the substantial codimension must be $p=n-1$ and, 
consequently, the normal bundle is flat. That $f$ is \emph{substantial} 
means that the codimension of $f$ cannot be reduced, that is, 
its image is not contained in any proper totally geodesic 
submanifold of the ambient space.\vspace{1ex}

\begin{theorem}\label{main} Let 
$f\colon M^n_c\to\Q_{\tilde c}^{n+p}$, $n\geq 2$ and $c\neq\tilde{c}$, 
be a substantial isometric immersion with $\dim N_1^f=n-1$ at any point.  
If $c>\tilde{c}$ assume also that $f$ is free of weak-umbilic
points. Then $p=n-1$.
\end{theorem}

We point out that for the class of  isometric immersions 
$f\colon M^n_c\to\Q^{2n-1}_{\tilde c}$ with $c\neq\tilde{c}$ and 
free of weak-umbilic points if $c>\tilde c$, there is an abundance of 
knowledge in the literature that makes this case of special 
interest. One basic property is that these submanifolds are at least 
locally holonomic. This means that $M^n_c$ carries local orthogonal 
coordinates such that the coordinate vector fields are 
principal directions. In particular, this is used to generalize
a classical correspondence for surfaces by showing that such isometric 
immersions are in correspondence with solutions of certain systems 
of partial differential equations. These are called the Generalized 
sine-Gordon, Generalized sinh-Gordon, Generalized Laplace, and Generalized 
wave equation, since these are generalizations of the classical 
two-dimensional versions; see Chapter $5$ of \cite{dt}. Moreover, 
given one such immersion, one can create a large family of new 
examples (sometimes even parametrically as shown in \cite{dt1}) 
by means of the B\"{a}cklund transformation (\cite{tete}, \cite{te}) 
or the Ribaucour transformation (\cite{dt1}).\vspace{1ex}

The following example shows that any substantial codimension is 
possible for submanifolds $f\colon M^n_c\to\Q^{n+p}_{\tilde c}$ 
with flat normal bundle if $N_1^f$ has the highest possible 
rank, namely, if $\dim N_1^f=n$. 

\begin{example} {\em Let $a_i$, $1\leq i\leq n$, be positive numbers 
such that $\sum_{i\geq 1}a_i^2=1$ and let $c\colon\R\to\Sf^m_{1/a_1^2}$, 
$m\geq 3$, be a spherical curve parametrized by arc length such that 
all the curvatures with respect to a Frenet frame never vanish. Then 
let $f\colon\R^n\to\Sf^{2n+m-2}_1\subset\R^{2n+m-1}$ be the isometric 
immersion given by
$$
f(t)=(c(t_1),a_2\sin (t_2/a_2),a_2\cos(t_2/a_2),\ldots, 
a_n\sin(t_n/a_n), a_n\cos(t_n/a_n))
$$
where $t=(t_1,\ldots,t_n)$. It is easy to verify  that $\dim N_1^f=n$ 
and that $f$ has flat normal bundle and is substantial.
}\end{example}

\section{The proof}

We first recall some basic facts about submanifolds with flat
normal bundle which can be seen in \cite{dt}.
\vspace{1ex}

If an isometric immersion $f\colon M^n\to\Q_c^{n+p}$ has flat 
normal bundle, that is, if the curvature tensor of the normal 
connection vanishes, it is a standard fact that at any point 
$x\in M^n$ there exists a set of unique pairwise distinct vectors 
$\eta_i(x)\in N_fM(x),\,1\leq i\leq s(x)$, called the 
\emph{principal normals} of $f$ at $x$ and an associate 
orthogonal splitting of the tangent space as
$$
T_xM=E_{\eta_1}(x)\oplus\cdots\oplus E_{\eta_s}(x),
$$
where 
$$
E_{\eta_i}(x)=\big\{X\in T_xM:\a_f(X,Y)=\<X,Y\>\eta_i\;\;
\text{for all}\;\;Y\in T_xM\big\}.
$$
Hence the second fundamental form of $f$ has the simple representation
\be\label{sfffnb}
\a_f(X,Y)=\sum_{i=1}^s\<X^i,Y^i\>\eta_i
\ee
where $X\mapsto X^i$ is the orthogonal projection onto $E_{\eta_i}$.

The dimension of $E_{\eta_i}(x)$ is called the \emph{multiplicity} of 
$\eta_i$  of $f$ at $x\in M^n$. If $s(x)=k$ is constant on $M^n$, the 
maps $x\in M^n\mapsto\eta_i(x)$, $1\leq i\leq k$, are smooth vector 
fields, called the \emph{principal normal vector fields} of $f$, and 
the distributions $x\in M^n\mapsto E_{\eta_i}(x)$, $1\leq i\leq k,$ 
are also smooth. \vspace{2ex}

\noindent\emph{Proof of Theorem \ref{main}:} We first observe that
the fact that $f$ has to have flat normal bundle is somehow inside
the arguments by Cartan and Moore. For a proof of this fact we refer 
to Theorem $5.5$ in \cite{dt}, where the statement is for codimension 
$n-1$. But it is straightforward to verify that the same conclusion 
holds under the weaker assumption that $\dim N_1^f=n-1$ since the 
proof reduces to analyze the algebraic structure of the second 
fundamental form as a map $\a_f\colon TM\times TM\to N_1^f$.

We show next that at any point there are exactly $n$ 
distinct principal normals $\eta_1,\dots, \eta_n$.  If $c<\tilde c$ 
this is trivial since if a principal normal $\eta_j$ has multiplicity 
$m>1$, we would have from the Gauss equation that 
$\|\eta_j\|^2=c-\tilde c<0$. Hence we may assume $c>\tilde c$. 
First observe that there is at most one principal 
normal of multiplicity $m>1$. In fact, if $\eta_1\neq\eta_2$
have both this property, using that
$\|\eta_1\|^2= c-\tilde c=\|\eta_2\|^2$ we would conclude 
that $\eta_1=\eta_2$.

Suppose that there is one principal normal $\eta_1$ of 
multiplicity $m>1$. Using \eqref{sfffnb} we have
$$
\<\a_f(X,Y),\eta_1\>
=\sum_{i=1}^{n-m+1}\<X^i,Y^i\>\<\eta_i,\eta_1\>
=(c-\tilde c)\<X,Y\>,
$$
and this is in contradiction with the assumption on weak-umbilics
points.

We claim that 
$$
(i)\;\dim N_1^f\geq n-1 \;\;\text{and}
\;\;(ii)\;\|\eta_i\|^2\neq c-\tilde c
\;\;\text{for all}\;\; 1\leq i\leq n.
$$ 
If $c<\tilde c$ then $(ii)$ is trivial and $(i)$ follows as an 
application of either Otsuki's lemma (see Corollary $6.2$ of \cite{dt}) 
or of the theory of flat bilinear forms (see Theorem $5.1$ of \cite{dt}). 
If $c>\tilde{c}$  suppose that
\be\label{eta}
\|\eta_1\|^2= c-\tilde c.
\ee
Then
\be\label{principal}
\|\eta_j\|^2\neq c-\tilde c\;\;\text{for all}\;\; 2\leq j\leq n
\ee
since if $\|\eta_j\|^2=c-\tilde c$ for some $2\leq j\leq n$, then
that $\<\eta_1,\eta_j\>=c-\tilde c$ easily gives $\eta_1=\eta_j$, 
a contradiction. Hence, we always have principal normals
$\eta_2,\dots,\eta_n$ such that \eqref{principal} holds. These
vectors  are linearly independent, which proves $(i)$. In fact, 
if $\sum_{i\geq 2}k_i\eta_i=0$, then
taking the inner product with $\eta_1$ gives 
$\sum_{i\geq 2} k_i=0$ whereas the inner product with $\eta_j$, 
$2\leq j\leq n$, yields 
$$
k_j(\|\eta_j\|^2+\tilde c-c)=(\tilde c -c)\sum_{i\geq 2}k_i=0.
$$
It follows using  \eqref{principal} that $k_j=0$ for all 
$2\leq j\leq n$. 

Since $\dim N_1^f=n-1$, set $\eta_1=\sum_{i\geq 2}k_i \eta_i$. 
Taking the inner product with $\eta_1$ and using \eqref{eta} 
gives
\be\label{kappa}
\sum_{i\geq 2}k_i=1.
\ee
Taking the inner product with $\eta_i,\ 2\leq i\leq n,$ and 
using \eqref{principal} and \eqref{kappa} we obtain 
$k_i=0$ in contradiction with \eqref{kappa}, and the proof 
of the claim is complete. \vspace{1ex}

Next we show that  the subbundle $N_1^f$ is parallel 
with respect to the normal connection. To see this, we
first observe that the Codazzi equation gives
\be\label{c}
\nabla_{X_j}^\perp\eta_i
=\<\nabla_{X_i} X_i,X_j\>(\eta_i-\eta_j),\; i\neq j
\ee
where the $X_i\in\Gamma(E_{\eta_i})$, $1\leq i\leq n$, are unit 
local vector fields.

Assume first that $n=2$. Then the principal normals 
$\eta_1,\eta_2$ are linearly dependent, in fact, they satisfy 
$$
(c-\tilde c)\eta_i-\|\eta_i\|^2\eta_j=0,\;\; i\neq j, 
$$ 
and the parallelism follows using \eqref{c}.

Next let $n\geq 3$. By the claim, there is a set of principal normals
$\eta_1,\dots,\eta_{n-1}$ that are linearly independent
and 
\be\label{eq3}
\eta_n=\sum_{i=1}^{n-1}\rho_i\eta_i.
\ee
Taking the inner product with $\eta_n$ gives
$$
\|\eta_n\|^2=(c-\tilde c)\sum_{i=1}^{n-1}\rho_i
$$
whereas with $\eta_j$, $1\leq j\leq n-1$, yields
$$
\rho_j\|\eta_j\|^2
=(\tilde c-c)\big(\sum_{i\neq j}\rho_i-1\big),\; 1\leq j\leq n-1.
$$
Then
$$
\rho_j=-\frac{\|\eta_n\|^2
+\tilde c- c}{\|\eta_j\|^2+\tilde c- c},\;1\leq j\leq n-1.
$$
It now follows from \eqref{eq3} that
\be\label{i}
\sum_{i=1}^n \frac{\eta_i}{\|\eta_i\|^2+\tilde c-c}=0.
\ee
Taking the inner product with some $\eta_j$ yields
\be\label{ii}
\sum_{i=1}^n \frac{1}{\|\eta_i\|^2+\tilde c-c}=\frac{1}{\tilde c-c}\cdot
\ee
It follows from \eqref{i} and \eqref{ii} that
\be\label{con}
\eta_i=\sum_{j\neq i}^n\frac{\tilde c-c}{\|\eta_j\|^2
+\tilde c-c}(\eta_i-\eta_j), \; 1\leq i\leq n.
\ee

We obtain using \eqref{con} and \eqref{c} that
$$
\left(1-\sum_{j\neq i}^n\frac{\tilde c-c}{\|\eta_j\|^2
+\tilde c-c}\right)\nabla_{X_i}^\perp\eta_i\in N_1^f.
$$
Since \eqref{ii} gives 
$$
1-\sum_{j\neq i}^n\frac{\tilde c-c}{\|\eta_j\|^2+\tilde c-c}
=\frac{\tilde c-c}{\|\eta_i\|^2+\tilde c-c}\neq 0,
$$
then $\nabla_{X_i}^\perp\eta_i\in N_1^f$, and using \eqref{c} 
it follows that $N_1^f$ is parallel in the normal connection. 

To conclude the proof, we recall an elementary fact from the 
theory of isometric immersions (cf.\ Proposition $2.1$ of 
\cite{dt}). If the first normal spaces form a parallel normal 
subbundle then the codimension reduces to the rank of $N_1^f$. 
Hence, in our case, we obtain that $p=n-1$.\qed

\noindent Marcos Dajczer\\
IMPA -- Estrada Dona Castorina, 110\\
22460--320, Rio de Janeiro -- Brazil\\
e-mail: marcos@impa.br
\bigskip

\noindent Christos-Raent Onti\\
Department of Mathematics and Statistics\\
University of Cyprus\\
1678, Nicosia -- Cyprus\\
e-mail: onti.christos-raent@ucy.ac.cy
\bigskip

\noindent Theodoros Vlachos\\
University of Ioannina \\
Department of Mathematics\\
Ioannina -- Greece\\
e-mail: tvlachos@uoi.gr

\end{document}